\documentclass{article}
\usepackage{amsmath, amssymb, amsthm}
\usepackage{geometry}
\usepackage{hyperref}
\usepackage[T1]{fontenc}

\usepackage{indentfirst}
\usepackage{graphicx} 

\begin{document}

\title{On the Composition of the Euler Function and the Dedekind Arithmetic Function}
\author{Aimin Guo, Huan Liu$^*$ and Qiyu Yang
\\ School of Mathematics and Statistics, Henan Normal University,
\\Xinxiang 453007, People's Republic of China.\\
\small{*Corresponding author : Huan Liu,
e-mail: liuhuan@htu.edu.cn}
}
\date{March 2025}
\maketitle

\begin{abstract}
Let $I(n) = \frac{\psi(\phi(n))}{\phi(\psi(n))}$ and $K(n) = \frac{\psi(\phi(n))}{\phi(\phi(n))}$, where $\phi(n)$ is Euler’s function and $\psi(n)$  is Dedekind’s arithmetic function. We obtain the maximal order of \(I(n)\), as well as the average orders of \(I(n)\) and \(K(n)\). Additionally, we prove a density theorem for both \(I(n)\) and \(K(n)\).
\end{abstract}

\maketitle

\section{Introduction}
 Let $\phi$ be Euler’s function and let $\psi$ be Dedekind’s arithmetic function. Here
\[\phi \left ( n \right ) =n\prod_{p|n}\left ( 1-\frac{1}{p}  \right ),\quad  \quad  \psi\left ( n \right )=n\prod_{p|n}\left ( 1+\frac{1}{p}  \right ),\]
where $p$ ranges over the prime factors of $n$. For some other studies on $\psi$ and $\phi$, see reference \cite{ref11}. In 2005, Sándor \cite{ref12} investigated the relation between $\psi(\phi(n))$ and  
 $\phi(\psi(n))$ under specific conditions. 
 
In this paper, motivated by the phenomenon that $\psi(\phi(n)) - \phi(\psi(n))$ can take both positive and negative values,
we investigate the expressions
\[I(n) := \frac{\psi(\phi(n))}{\phi(\psi(n))}\quad \text{and} \quad 
K(n) := \frac{\psi(\phi(n))}{\phi(\phi(n))},\]
to further reveal the intrinsic relationship between $\psi(\phi(n))$ and $\phi(\psi(n))$.

In particular, we establish three key results: first, the maximal order of \( I(n) \); second, the average orders of both \( I(n) \) and \( K(n) \); and third, the density theorems for \( I(n) \) and \( K(n) \). Notably, from the density theorem for \( I(n) \), we further demonstrate that the arithmetic difference $\psi(\phi(n)) - \phi(\psi(n))$ is predominantly positive for most positive integers $n$.

If $ k $ is a positive integer, we denote by $ \log_k x $ the $ k $-th iteration of the logarithm. Throughout this paper, the letters $ p $, $ q $, and $ r $ represent prime numbers; $ \gamma $ denotes the Euler–Mascheroni constant; $ \varepsilon $ is any  positive parameter that may not the same; and $ \omega(n) $ represents the number of distinct prime factors of $ n $.

Based on our investigation, we establish the following theorems:

\par\noindent\textbf{Theorem 1.1} 
The maximal order of $ I(n) $ is characterized by
    \begin{equation}
        \limsup_{n \to \infty} \frac{I(n)}{\log_2^2 n} = \frac{6}{\pi^2} e^{2\gamma}. \label{eq:maximal}\tag{1.1}
    \end{equation}
    
\par\noindent\textbf{Theorem 1.2}

\begin{enumerate}
    \renewcommand{\labelenumi}{(\arabic{enumi})} 
    \item 
    As $x \to \infty$, we have
   \begin{equation}
   \frac{1}{x} \sum_{n \leq x} I(n) = c_0\frac{6}{\pi^2}e^{2\gamma} \log_{3}^{2} x + O(\log_{3}^{3/2} x),\label{300}\tag{1.2}
  \end{equation}
  \hspace*{-2.5em} 
    \noindent
where 

\[
c_0 =\lim _{x \rightarrow \infty} \frac{1}{x} \sum_{n \leq x} \frac{\phi(n)}{\psi(n)}= \prod_{p} \left(1 - \frac{2}{p(p+1)}\right)\approx 0.47.\]
  \item It follows that on a set of density 1,
 \begin{equation}
 I(n)= \left(1 + o\left(1\right)\right) \frac{6}{\pi^2}  e^{2\gamma} \log_3^2 n  \prod_{p|n} \frac{p-1}{p+1}.
  \label{301}\tag{1.3}
\end{equation}
 
\end{enumerate}

In particular, \eqref{301} shows that $\psi(\phi(n)) - \phi(\psi(n))$ is predominantly positive for most positive integers $n$.

 Following a similar methodology applied to $I(n)$, we investigate $K(n)$, establish its average order, and prove the density theorem.

\par\noindent\textbf{Theorem 1.3} 
\begin{enumerate}
    \renewcommand{\labelenumi}{(\arabic{enumi})} 
    
    \item As \( x \to \infty \), the average order of \( K(n) \) is
    \begin{equation}
        \frac{1}{x} \sum_{n \leq x} K(n) = \frac{6}{\pi^2} e^{2\gamma} \log_3^2 x + O(\log_3^{3/2} x). \label{201}\tag{1.4}
    \end{equation}
    
      \item It follows that on a set of density 1,
      \begin{equation}
     K(n)= \left(1 + o\left(1\right)\right) \frac{6}{\pi^2}  e^{2\gamma} \log_3^2 n 
     .\label{302}\tag{1.5}
\end{equation}
\end{enumerate}

\section{Preliminary Results}
For all sufficiently large $x$, we present the following known standard results in number theory:
\begin{align}
    \prod_{p \leq x} \left(1 - \frac{1}{p}\right) 
    &= \frac{e^{-\gamma}}{\log x} \left(1 + O\left(\frac{1}{\log x}\right)\right),
    \label{3-1}
    \tag{2.1} \\
    \prod_{p \leq x} \left(1 + \frac{1}{p}\right) 
    &= \frac{6e^{\gamma}}{\pi^2} \log x \left(1 + O\left(\frac{1}{\log x}\right)\right).
    \label{3-2}
    \tag{2.2}
\end{align}

\textbf{Proof.} Equation (\ref{3-1}) can be found in Mertens’ work \cite{ref9}, and Equation (\ref{3-2}) was derived from Littlewood’s paper, as demonstrated in \cite{ref7}.

\textbf{Lemma 2.1} 
\[
\liminf_{n \to \infty} \frac{\phi(n) \log_2 n}{n} = e^{-\gamma},\label{3-3}\tag{2.3}
\]
\[\limsup_{n \to \infty} \frac{\psi(\phi(n))}{n \log_2 n} = \frac{6}{\pi^2} e^\gamma,\label{3-4}\tag{2.4}
\]
\[
\limsup_{n \to \infty} \frac{\psi(n)}{n \log_2 n} = \frac{6}{\pi^2} e^\gamma,\label{3-5}\tag{2.5}
\]

\[\liminf_{n \to \infty} \frac{\phi(\psi(n)) \log_2 n}{n} = e^{-\gamma},\label{3-6}\tag{2.6}\]

\[
\limsup_{n \to \infty} \frac{\psi(\phi(n))}{\phi(\phi(n))\log_2^2 n} = \frac{6}{\pi^2} e^{2\gamma}.
\label{3-7}\tag{2.7}\]

\textbf{Proof.}
For \eqref{3-3} we refer to Landau \cite{ref6}, while Eqs.~\eqref{3-4}--\eqref{3-7} can be found in Sándor \cite{ref10}.

\textbf{Lemma 2.2} 
Suppose $ k $ and $ a $ are coprime positive integers. There exists a prime $ p $ such that
$p \equiv a \pmod{k}$ and $p = O(k^{11/2}).
$

\textbf{Proof.}
It follows from Theorem 6 of Heath-Brown~\cite{ref4}.

\textbf{Lemma 2.3} 
Suppose that $x$ is sufficiently large. Then there exists a constant $c_1$ such that, for all positive integers $n<x$ with $O\left(\frac{x}{\log _{3}^2 x}\right)$  exceptions, both  $\psi(n)$  and  $\phi(n)$  are divisible by all prime powers  $p^{a}$  satisfying  $p^{a}<c_{1}\frac{\log _{2} x}{\log _{3} x}$.

\textbf{Proof.}
The above result for the case of the function $\phi(n)$ is Lemma 2 in \cite{ref8}. 
To prove this Lemma, we consider the number of $ n $ for which $ p^a \nmid \psi(n) $. 

Since $\psi(n)$ is a multiplicative function, if $p^a \nmid \psi(n)$, then there exists at least one prime factor $q$ such that $p^a \nmid \psi(q^e)$. 
Therefore, we only need to consider the case for $p^a\nmid q+1$.

To prove the desired result, we employ the method from Lemma 4 in \cite{ref2}. Consider an arbitrary positive integer \( m \) and a positive real parameter \( x \). Define 
\[
S(x, m) = \sum_{\substack{q \leq x \\ m \,|\, (q+1)}} \frac{1}{q}.
\]
Through application of the Siegel-Walfisz theorem combined with partial summation techniques, we derive that there exist absolute constants \(c_1 > 0\) and \(x_0 > 0\) such that for all \(x > x_0\) and integers \(m \leq \log x\), 

\[
S(x, m) \geq \frac{c_1 \log_2 x}{\phi(m)}.
\]
We define the exceptional set \(\mathcal{N}_m\) as the family of integers \(n \leq x\) that satisfy the following condition: there exists no prime factor \(q\) of \(n\) with the congruence relation
\[
q \equiv -1 \pmod{m}.
\]
Next we apply Brun's sieve. Assume \( x_0 \) is chosen sufficiently large such that for all \( x > x_0 \), the inequality 
$\log x > g(x) = c_1 \frac{\log_2 x}{\log_3 x}$
holds. Taking \( m = p^a < g(x) \), we derive that
\[\# \mathcal{N}_{p^{a}}=\sum_{\substack{n \leq x \\ p^{a} \nmid q+1}} 1<\frac{c_{2} x}{\exp \left(S\left(x, p^{a}\right)\right)}<\frac{c_{2} x}{\exp \left(\log _{3} x\right)}=\frac{c_{2} x}{\log _{2} x} .\]
Here $c_2$ is a positive constant.  
Therefore
\[\sum_{p^{a}<g(x)} \sum_{\substack{n \leq x \\ p^{a} \nmid q+1}} 1\ll\frac{x g(x)}{\log_{2} x \cdot \log g(x)}\ll \frac{x}{\log_{3}^{2} x}.\]
Thus, we obtain that among all the integers $n<x$, the number of interges $n$ such that $\psi(n)$ is not divisible by all the prime powers $ p^a < g(x)$ is at most $\frac{x}{\log_{3}^{2} x}$. This completes the proof of Lemma 2.3.

\textbf{Lemma 2.4}
As $x \to \infty$,
\[
\sum_{n \leq x} \frac{\phi(n)}{\psi(n)} = c_0 x + O(x^{\varepsilon}),
\]
where $ c_0 $ is the constant from Theorem~1.2 and $ \mu(n) $ denotes the Möbius function.

\textbf{Proof.} Let $\zeta(s)$ denote the Riemann zeta function. Then we have
\begin{equation}
\begin{aligned}
\sum_{n=1}^\infty \frac{\phi(n)/\psi(n)}{n^s} &=
\prod_p \left( 1 + \frac{\frac{p - 1}{p+1}}{p^{s}} + \frac{\frac{p - 1}{p+1}}{p^{2s} } + \frac{\frac{p - 1}{p+1}}{p^{3s}} + \cdots \right)\\
&= \zeta(s) \prod_p \left( 1 - \frac{1}{p^s} \right) \prod_p \left( 1 + \frac{\frac{p - 1}{p+1}}{p^{s}} + \frac{\frac{p - 1}{p+1}}{p^{2s} } + \frac{\frac{p - 1}{p+1}}{p^{3s}} + \cdots \right)\\
&=\zeta(s)\prod_p \left( 1 + \frac{\frac{p-1}{p+1}-1}{p^s} +\frac{\frac{p - 1}{p+1}-\frac{p - 1}{p+1}}{p^{2s} }+ \cdots \right)\\
&=\zeta(s)\prod_p \left( 1 + \frac{\frac{p-1}{p+1}-1}{p^s}  \right)\\
&= \zeta(s) R(s),\nonumber
\end{aligned}
\end{equation}
where Dirichlet series  \( R(s) \) converges absolutely in the half-plane \( \Re(s)>0 \). Write
\[
R(s) = \sum_{n=1}^{\infty} \frac{a_n}{n^s}, a_n=\mu ^{2}\left ( n \right ) \prod_{p|n}\frac{-2}{p+1},
\]
then $\frac{\phi(n)}{\psi(n)}=\sum_{d|n}a_d.$ Therefore
\begin{equation}
\begin{aligned}
\sum_{n \leq x} \frac{\phi(n)}{\psi(n)} & =\sum_{n \leq x} \sum_{d \mid n} a_{d}=\sum_{n \leq x} a_{d}\left[\frac{x}{d}\right]=x \sum_{d \leq x} \frac{a_{d}}{d}+O\left(\sum_{d \leq x}\left|a_{d}\right|\right) \\
& =R(1) x+O\left(x \sum_{d>x} \frac{\left|a_{d}\right|}{d}\right)+O\left(\sum_{d \leq x}\left|a_{d}\right|\right)\nonumber,
\end{aligned}
\end{equation}
where \( R(1) = c_0 \).

Since
\[
\sum_{d \le x} |a_d| = \sum_{d \le x} \frac{|a_d|}{d^{\varepsilon }} d^{\varepsilon } = O(x^{\varepsilon })
\]
and
\[
\sum_{d > x} \frac{|a_d|}{d} 
= \sum_{d > x} \frac{|a_d|}{d^{\varepsilon } d^{1-\varepsilon }} 
\le x^{-1+\varepsilon } \sum_{d > x} \frac{|a_d|}{d^{\varepsilon }} = O(x^{-1+\varepsilon }),
\]
it follows that
\[
\sum_{n \le x} \frac{\phi(n)}{\psi(n)} = R(1)x + O(x^{\varepsilon}),
\]
which completes the proof of Lemma 2.4.

\textbf{Lemma 2.5}
 Let $x$ be a positive real number. Setting
\[  h_\phi\left(n\right)=\sum_{\substack{p| \phi \left ( n \right )  \\ p> \log_{2}x  }}\frac{1}{p} \quad \text{and} \quad h_\psi\left(n\right)=\sum_{\substack{p| \psi \left ( n \right )  \\ p> \log_{2}x  }}\frac{1}{p},\]
then
\[\sum_{n\le x} h_\phi \left(n\right)\ll \frac{x}{\log_{3}x } \quad \text{and} \quad \sum_{n\le x} h_\psi\left(n\right)\ll \frac{x}{\log_{3}x }.\]
The inequalities stated below are valid on a set of asymptotic density 1. 
\[h_\psi\left(n\right)=\sum_{\substack{p| \psi \left ( n \right )  \\ p> \log_{2}x  }}\frac{1}{p}< \frac{\log_4(n)}{\log_3(n)}  \quad \text{and} \quad  h_\phi\left(n\right)=\sum_{\substack{p| \phi \left ( n \right )  \\ p> \log_{2}x  }}\frac{1}{p}< \frac{\log_4(n)}{\log_3(n)}.\]

\textbf{Proof.}
 Clearly we have, for a fixed prime
\[\sum_{\substack{n\leq x\\ p| \psi \left ( n \right ) }}1=\frac{x}{p^2} +\sum_{\substack{q\leq x\\ p| q+1} }\frac{x}{q }\ll \frac{x}{p^2} +\frac{x\log_{2} x}{p}\ll \frac{x\log_{2} x}{p}.\]
It now follows that
\begin{equation}
\sum_{n \leq x} h_\psi(n) = \sum_{p \leq x} \frac{1}{p} \sum_{\substack{n\leq x\\ p| \psi \left ( n \right ) }} 1
\ll x \log_2 x \sum_{p > \log_2 x} \frac{1}{p^2}
\ll \frac{x}{\log_3 x}.\label{eq:22}\nonumber
\end{equation}
Similarly
\begin{equation}
\sum_{n \leq x} h_\phi(n) 
\ll \frac{x}{\log_3 x}.\label{eq:23}\nonumber
\end{equation}

We define $ x_n = h_\psi(n) $ and $ \mathbb{E}[x_n] = \frac{1}{x} \sum_{n \leq x} h_\psi(n)= O\left( \frac{1}{\log_3 x} \right) $. Applying Markov's inequality, we obtain that for all $ \varepsilon > 0 $,
$$
P\left( x_n \ge \varepsilon \right) \le \frac{\mathbb{E}[x_n]}{\varepsilon}.
$$
Choosing $ \varepsilon = \frac{\log_4 x}{\log_3 x} $, we derive that
$$
P\left( h_\psi(n) \ge \frac{\log_4 x}{\log_3 x} \right) \le \frac{O\left( \frac{1}{\log_3 x} \right)}{\frac{\log_4 x}{\log_3 x}} = O\left( \frac{1}{\log_4 x} \right).
$$

So
\[\# \left\{ n \leq x : h_\psi(n) \geq \frac{\log_4 x}{\log_3 x} \right\} \leq O \left( \frac{x}{\log_4 x} \right) = o(x).\]
The proportion of integers $ n \leq x $ satisfying $ h_\psi(n) < \frac{\log_4 x}{\log_3 x} $ is  
$$
1 - \frac{o(x)}{x} \to 1 \quad (x \to \infty).
$$

Therefore, the inequality $ h_\psi(n) < \frac{\log_4 x}{\log_3 x} $ holds on a set of asymptotic density 1.

Similarly, we obtain $h_\phi(n) < \frac{\log_4 (n)}{\log_3 (n)}$  holds on a set of asymptotic density 1, the proof of Lemma 2.5 is complete.

\textbf{Lemma 2.6}
There exists a constant $b_3$ such that the set of positive integers 
$n \le x$ for which $\omega(\psi(n)) > b_3 \log_2^2 x$ contains at most 
$O\bigl(\frac{x}{\log x}\bigr)$ elements. The conclusion remains valid when $\psi(n)$ is replaced by $\phi(n)$.

\textbf{Proof.}
In the paper by Jean-Marie De Koninck and Florian Luca \cite{ref2}, the result for $\phi$ has already been provided. We will adopt the same approach used by them to prove $\psi$.

We begin by noting that, according to the definition of $\omega \left ( \psi \left ( n \right )  \right )$, to ensure $\omega \left ( \psi \left ( n \right )  \right )$ is sufficiently large, it suffices to constrain both the number and the size of prime factors in $n$. We now proceed to classify 
$n$ as follows:

 First let $\mathcal{D}_1$ be the set of all $n \leq x$ such that $k = \omega(n) > 3e \log_2 x$. The Hardy--Ramanujan theorem, established in the seminal work \cite{ref5}, is a fundamental result in analytic number theory and states that
\[
\#\{n \leq x : \omega(n) = k\} \ll \frac{x}{\log x} \cdot \frac{1}{(k-1)!} \cdot (\log_2 x + O(1))^{k-1}.
\]
Combining Stirling's formula with the inequality yields the asymptotic relation:
\[
\#\{n \leq x : \omega(n) = k\} \ll \frac{x}{\log x} \cdot 
\left( \frac{e \log_2 x + O(1)}{k-1} \right)^{k-1} 
< \frac{x}{\log x} \cdot \frac{1}{2^{k-1}}.
\]
Since $k - 1 > 3e \log \log x - 1$ and $x$ is assumed to be large, 

\[
\#\mathcal{D}_1 = \#\{n \leq x : \omega(n) > 3e \log_2 x\} 
\ll \frac{x}{\log x} \sum_k \frac{1}{2^k} 
\ll \frac{x}{\log x} 
.
\]

Let $\mathcal{D}_2$ be the set of those $n \leq x$ which are divisible by a prime number $p$ such that $\omega(p-1) \geq b := \lfloor e^2 \log_2 x \rfloor$. Then
\[
\#\mathcal{D}_2 \leq \sum_{\substack{p \leq x \\ \omega(p-1) \geq b}} \frac{x}{p} 
\leq x \sum_{k \geq b} \frac{1}{k!} \left(\sum_{q^a \leq x} \frac{1}{q^a}\right)^k.
\]

Using the inequality
\[
\sum_{q^a \leq x} \frac{1}{q^a} \ll e \log_2 x + O(1),\]
we have
\[
\#\mathcal{D}_2 \ll x \sum_{k \geq b} \left(\frac{e \log_2 x + O(1)}{k}\right)^k 
\ll x \sum_{k \geq b} \frac{1}{2^k} 
\ll \frac{x}{\log x} 
.
\]
Finally, let $\mathcal{D}_3$ be the set of those $n \leq x$ which are divisible by a prime number $p$ such that $\omega(p+1) \geq b$. The same argument as above shows that
\[
\#\mathcal{D}_3 \ll \frac{x}{\log x}.
\]
If $n$ does not belong to the partition $D=\mathcal{D}_1\cup\mathcal{D}_2\cup\mathcal{D}_3 $
defined above, then it satisfies the following conditions:

1. The number of distinct prime factors is moderate:
$\omega\left ( n \right ) \le 3e\log_2 x $.

2. For all primes \( p \) dividing \( n \), the number of distinct prime factors of \( p \pm 1 \) satisfies
 $\omega(p \pm 1) \leq b = \left\lfloor e^2 \log_2 x \right\rfloor.$
 
 In summary, we conclude that when $n\notin D$
\[\omega \left ( \psi  \left ( n \right )  \right ) \le \omega \left ( n \right ) +\sum_{p|n}\omega \left (p+1  \right ) \le 3e\log_2 x+b\cdot 3e\log_2 x =O\left ( \log_2^{2}x \right ).\]

Therefore, the condition
$\omega(\psi(n)) > b_3 \log_2^2 x$ holds only if $n \in D$. In addition, $D$ contains $O(x/\log x)$ elements, as claimed.

\section{{Proof of Theorem 1.1}}

We will show that for $n$ sufficiently large,
\begin{equation}
I(n) \leq (1 + o(1)) \frac{6}{\pi^2} e^{2\gamma} \log_2^2 n.\label{eq:1}\tag{3.1}
\end{equation}
The proof of inequality (\ref{eq:maximal}) in Theorem 1.1 will therefore follow if there exists an infinite sequence $\{n_k\}$ with $n_k\to \infty$ such that
$I(n_k) \geq (1 + o(1)) \frac{6}{\pi^2} e^{2\gamma} \log_2^2 n_k.$

To prove (\ref{eq:1}), we first derive from (\ref{3-4}) that
\begin{equation}
\psi(\phi(n)) \leq (1 + o(1))\frac{6}{\pi^2} e^\gamma n \log_2 n.\label{eq:2}\tag{3.2}
\end{equation}
Additionally, it can be obtained from equation (\ref{3-6}) that 
\begin{equation}
\phi(\psi(n)) \geq (1 + o(1)) e^{-\gamma} \frac{n}{\log_2 n}.\label{eq:3}\tag{3.3}
\end{equation}
By combining (\ref{eq:2}) and (\ref{eq:3}), we obtain (\ref{eq:1}).

 To complete the proof of Theorem 1.1, we construct an infinite sequence $\{n_k\}$ with $n_k \to \infty$ such that
$I(n_k) \geq (1 + o(1)) \frac{6}{\pi^2} e^{2\gamma} \log_2^2 n_k.$ So let $x$ be a large integer, and let $P$ and $Q$ be the smallest
primes such that
\[P \equiv 1 \pmod{M(x)} \quad \text{and} \quad Q \equiv -1 \pmod{M(x)},
\]
where $M(x) = \mathrm{LCM}[1, 2, \dots, x]$, and set
\[n = P Q.\]
By employing the prime number theorem, we can derive that
\begin{equation}
\begin{aligned}
M(x) =\prod_{p^{k}\le x}p^{k}=\exp({\sum_{p^{k}\le x}  \log p^{k}})&=\exp{\left ( \sum_{p\le x} \log p+\frac{\log x}{\log 2}\sum_{p\le \sqrt{x} } \log p  \right ) }\\
= e^{(1 + o(1)) x} < e^{2x}.\nonumber
\end{aligned}
\end{equation}
Furthermore, by applying Lemma 2.2, we obtain
\[P \ll e^{11 x}, \quad Q \ll e^{11 x}.\]
Therefore, when $x$ is sufficiently large, we have the inequality $n < e^{23x}$. For our selection of a specially chosen positive integer $n$, we have,
since $\phi(n) = (P - 1)(Q - 1)$,
\begin{equation}\label{eq:4} \tag{3.4}
\begin{aligned}
      \frac{\psi(\phi(n))}{\phi(n)} &=\prod_{p\mid \phi(n)} \left( 1 + \frac{1}{p} \right)=\prod_{p\mid (P-1)(Q-1)} \left( 1 + \frac{1}{p} \right)\\
     &\geq \prod_{p\mid (P-1)} \left( 1 + \frac{1}{p}  \right)\\
     &\geq \prod_{p \mid M(x)} \left( 1 + \frac{1}{p} \right)
\geq\prod_{p \leq x}\left( 1 + \frac{1}{p}  \right).
\end{aligned}
\end{equation}
Using this in (\ref{eq:4}), we obtain that, by (\ref{3-2}),

\begin{equation}
\frac{\psi(\phi(n))}{\phi(n)} 
\geq \frac{6e^\gamma}{\pi^2}  \log x  \left( 1 + o\left({1}\right) \right).\label{eq:5}\tag{3.5}
\end{equation}
On the other hand, we know $\psi(n) = (P + 1)(Q + 1)$, so that
\begin{equation}\label{eq:6}\tag{3.6}
\begin{aligned}
\frac{\phi(\psi(n))}{\psi(n)} &= \prod_{p \mid (P+1)(Q+1)} \left( 1 - \frac{1}{p} \right) \leq \prod_{p \mid (P+1)} \left( 1 - \frac{1}{p} \right)\\ 
&\leq \prod_{p \mid M(x)} \left( 1 - \frac{1}{p} \right)
= \prod_{p \leq x} \left( 1 - \frac{1}{p} \right) 
= (1 + o(1)) \frac{e^{-\gamma}}{\log x},
\end{aligned}
\end{equation}
where we used (\ref{3-1}).

Gathering (\ref{eq:5}) and (\ref{eq:6}), we get that
\begin{equation}
I(n) \cdot \frac{\psi(n)}{\phi(n)} \geq (1 + o(1))\frac{6e^\gamma}{\pi^2}  e^{2\gamma} \log^2 x.
\label{eq:7}\tag{3.7}
\end{equation}

For chosen $n$, we obtain 
\begin{equation}\exp\{(1 + o(1))x\} < n < \exp\{23x\}.\label{eq:500}\tag{3.8}
\end{equation}
By taking the logarithm of both sides of equation (\ref{eq:500}), we obtain $\log_2 n = \log x + O(1)$. Therefore, equation \eqref{eq:7} can be written as
\begin{equation}
I(n) \cdot \frac{\psi(n)}{\phi(n)} \geq (1 + o(1)) \frac{6e^\gamma}{\pi^2} e^{2\gamma} \log_2^2 n.
\nonumber
\end{equation}

As $ x \to \infty $, since $ M(x) $ divides both $ P-1 $ and $ Q+1 $, it follows that $ P $ and $ Q $ also tend to infinity. We derive
\[\frac{\psi(n)}{\phi(n)} = \frac{(P + 1)(Q + 1)}{(P - 1)(Q - 1)} = 1 + o(1).\]
By combining equation (\ref{eq:1}) with the sequence we constructed, where the value of $n$ varies with $x$, we finally complete the proof of Theorem 1.1.

\section{Proof of Theorem 1.2}
\subsection{Proof of (\ref{300})}
We extend the approach in \cite{ref8} to study the mean value of $I(n)$. Define $M_0(x)=\mathrm{LCM}\left\{p^{a} \mid p^{a}<g(x)\right\}$, where $g(x) = c_1 \log_2 x / \log_3 x$ 
and the constant $c_1$ is specified in Lemma 2.3. Moreover, let 
\[
\mathcal{A}=\mathcal{A}(x) = \{ n : \sqrt{x} < n \leq x \text{ and } M_0(n) \mid \gcd(\phi(n), \psi(n)) \}.
\]

Combining (\ref{3-2}), we obtain
\begin{equation}
\frac{\psi(\phi(n))}{\phi(n)}\geq\frac{\psi(M_0(n))}{M_0(n)} = 
\prod_{p \leq g(n)} \left(1 + \frac{1}{p} \right)=\frac{6e^\gamma}{\pi^2} \log_3 x + O\left(\log_4 x\right) 
 \quad (n \in \mathcal{A}).
\label{eq:16}\tag{4.1}
\end{equation}

By following the analogous approach and subsequently applying (\ref{3-1}), we derive 
\begin{equation}
\frac{\phi(\psi(n))}{\psi(n)} \leq \frac{\phi(M_0(n))}{M_0(n)} = \prod_{p \leq g(n)} \left(1 - \frac{1}{p}\right) = \frac{e^{-\gamma}}{\log_3 x} \left(1 + O\left(\frac{\log_4 x}{\log_3 x}\right)\right) \quad (n \in \mathcal{A}).
\label{eq:17}\tag{4.2}
\end{equation}

Combining equations \eqref{eq:16} and \eqref{eq:17}, we derive the following result for $n \in \mathcal{A}$:

\begin{equation}
I(n) \geq  \frac{\phi(n)}{\psi(n)}\left(\frac{6}{\pi^2}  e^{2\gamma} \log_3^2 x 
 + O\left( \log_3 x\cdot\log_4 x \right)\right).
\label{eq:18}\tag{4.3}
\end{equation}
Summing over equation \eqref{eq:18}, we obtain

\begin{equation}
\sum_{n \leq x} I(n) \geq 
\sum_{\substack{n \in \mathcal{A} \\ n \leq x}}I(n) \geq \left(\frac{6}{\pi^2}e^{2\gamma} 
\log_3^2 x  + O\left( \log_3 x\cdot\log_4 x   \right) \right)\sum_{\substack{n \in \mathcal{A} \\ n \leq x}} \frac{\phi(n)}{\psi(n)}.
\label{eq:19}\tag{4.4}
\end{equation}
Next, we apply Lemma~2.3 to estimate the cardinality of the set $\{n \in [1, x] \setminus A\}$, and then combine it with the inequality $\phi(n) \leq \psi(n)$ that holds for all $n$, to derive:
\[
 \sum_{\substack{n \in \mathcal{A} \\ n \leq x}} \frac{\phi(n)}{\psi(n)} \geq \sum_{n \leq x} \frac{\phi(n)}{\psi(n)} - (x - \#\mathcal{A}) = \sum_{n \leq x} \frac{\phi(n)}{\psi(n)} + O\left( \frac{x}{\log_3^2 x} \right)= c_0 x + O\left( \frac{x}{\log_3^2 x} \right),\]
where $c_0$ is the constant defined in Lemma 2.4. By combining \eqref{eq:19}, we derive
\begin{equation}
\sum_{n \leq x} I(n) \geq c_{0} \frac{6}{\pi^2}e^{2 \gamma} x \log _{3}^{2} x+O\left(x \log _{3} x \cdot \log_4 x\right) .
\label{eq:20}\tag{4.5}
\end{equation}

To determine the sum $\sum_{n \leq x} I(n)$, further analysis of the upper bound for the sum $\sum_{n \leq x} I(n)$ is required.

Motivated by \cite{ref2}, we employ the interval partitioning method to investigate the upper bound of $\sum_{n \leq x} I(n)$. First, for $n \in[1,\sqrt{x} ]$, by applying inequality (\ref{eq:maximal}), we establish that
\begin{equation}
\sum_{n \leq \sqrt{x}} I(n)=O\left(\sqrt{x} \log _{2}^{2} x\right).
\label{eq:21}\tag{4.6}
\end{equation}

Next, we consider the set
\[\mathcal{B}=\mathcal{B}(x)=\left\{n: \sqrt{x}<n \leq x, h_{\phi}(n)<\frac{1}{\sqrt{\log _{3} x}}, h_{\psi}(n)<\frac{1}{\sqrt{\log _{3} x}}\right\},\]
and for each integer  $n \in \mathcal{B},$ decompose $\phi(n)$  as $\phi(n)=n_{1} \cdot n_{2},$ where
\[n_{1}=\prod_{\substack{p^{\alpha_{p} }\|\phi(n) \\ p \leq \log _{2} x}} p^{\alpha_{p}} \quad \text { and } \quad n_{2}=\prod_{\substack{p^{\alpha_{p} }\|\phi(n) \\ p>\log _{2} x}} p^{\alpha_{p}}.\] 
So that, using (\ref{3-2}),
\begin{equation}\label{eq:26}\tag{4.7}
\begin{aligned}
\frac{\psi(\phi(n))}{\phi(n)} & =\prod_{p \mid n_{1}}\left(1+\frac{1}{p}\right) \cdot \prod_{p \mid n_{2}}\left(1+\frac{1}{p}\right) \\
& \leq\left(\frac{6}{\pi^2} e^{\gamma} \log _{3} x+O(1)\right) \cdot \exp \left(O\left(h_{\phi}(n)\right)\right) \\
& =\left(\frac{6}{\pi^2}e^{\gamma} \log _{3} x+O(1)\right) \cdot \exp \left\{O\left(\frac{1}{\sqrt{\log _{3} x}}\right)\right\} \\
& =\frac{6}{\pi^2}e^{\gamma} \log _{3} x+O\left(\sqrt{\log _{3} x}\right) \quad(n \in \mathcal{B}) .
\end{aligned}
\end{equation}
On the other hand, given $n \in \mathcal{B}$, decompose $\psi(n)$ as $\psi(n)=m_{1} \cdot m_{2} $, where
\[m_{1}=\prod_{\substack{p^{\alpha_{p} }\|\psi(n) \\ p \leq \log _{2} x}} p^{\alpha_{p}} \quad \text { and } \quad m_{2}=\prod_{\substack{p^{\alpha_{p} }\|\psi(n) \\ p>\log _{2} x}} p^{\alpha_{p}}.\] 
Using an analogous approach, we derive that

\begin{equation}\label{eq:27}\tag{4.8}
\begin{aligned}
\frac{\phi(\psi(n))}{\psi(n)} & =\frac{\phi(m_1)}{m_1} \cdot \frac{\phi(m_2)}{m_2}\\
&\ge \frac{e^{-\gamma}}{\log_3 x}\left(1 + O\left(\frac{1}{\log_3 x}\right)\right)
\left(1 + O\left(\frac{1}{\sqrt{\log_3 x}}\right)\right)
 \\
 &= \frac{e^{-\gamma}}{\log_3 x} 
\left(1 + O\left(\frac{1}{\sqrt{\log_3 x}}\right)\right)
\quad (n \in \mathcal{B}).
\end{aligned}
\end{equation}
Combing \eqref{eq:26} and \eqref{eq:27}, we obtain

\begin{equation}
I(n) \le \frac{\phi(n)}{\psi(n)} \frac{6}{\pi^2}e^{2\gamma} \log_3^2 x
\left(1 + O\left(\frac{1}{\sqrt{\log_3 x}}\right)\right)
\quad (n \in \mathcal{B}).\label{eq:28}\tag{4.9}
\end{equation}
It conclude that
\begin{equation}\label{eq:29}\tag{4.10}
\begin{aligned}
\sum_{\substack{n \le x \\ n \in \mathcal{B}}} I(n) 
&\le \frac{6}{\pi^2}e^{2\gamma} \log_3^2 x 
\left(1 + O\left(\frac{1}{\sqrt{\log_3 x}}\right)\right) 
\sum_{\substack{n \le x \\ n \in \mathcal{B}}} \frac{\phi(n)}{\psi(n)}\\
&\le \frac{6}{\pi^2}e^{2\gamma} c_0 \, x \log_3^2 x + O\bigl(x \log_3^{3/2} x\bigr).
\end{aligned}
\end{equation}

We proceed by partitioning the interval to analyze the contribution of integers 
$n \in [\sqrt{x}, x]$ not belonging to the set $\mathcal{B}$. Specifically, such integers must be contained within the union $\mathcal{C}_\phi \cup \mathcal{C}_\psi$.

Here, for each function$f \in \{\phi, \psi\}$, the set $\mathcal{C}_f$ is defined as the collection of all integers $n \in [\sqrt{x}, x]$ satisfying $h_f(n) \ge 1/\sqrt{\log_3 x}$.
 
From Lemma~2.5, we obtain
\[
\frac{x}{\log_3 x} \gg \sum_{n \in \mathcal{C}_f} h_f(n) 
\ge \frac{\#\mathcal{C}_f}{\sqrt{\log_3 x}},
\]
so that
\begin{equation}\label{eq:33}\tag{4.11}
\#\mathcal{C}_f \ll \frac{x}{\sqrt{\log_3 x}} 
\quad \text{for } f = \phi \text{ and } f = \psi.
\end{equation}
We split the positive integers in the interval $\mathcal{C}_\phi \cup \mathcal{C}_\psi$,  into
the following two subsets:
\begin{equation}\label{eq:34}\nonumber
\begin{aligned}
&E = \left\{ n \;\middle|\;n\in \mathcal{C}_\phi \cup \mathcal{C}_\psi ; n\in \mathcal{D}  \right\},\\
&F = \left\{ n \;\middle|\;n\in \mathcal{C}_\phi \cup \mathcal{C}_\psi ; n\notin  \mathcal{D}  \right\},
\end{aligned}
\end{equation}
where 
$D$ is the set defined in Lemma 2.6.

Next, we only need to estimate the sum of $I(n)$ over the sets $E$ and $F$.
 Since by (\ref{eq:maximal}), \( I(n) \ll \log_2^2{n}\ll \log_2^2{x} \), it follows that
\begin{equation}\label{eq:35}\tag{4.12}
\sum_{n \in E} I(n) = O\left (\frac{x\log_2^2 x}{\log x}\right ).
\end{equation}
 It is straightforward to verify that if \( n \in F \), then both inequalities \( \omega(\phi(n)) < b_3 \log_2^2{x} \) and \( \omega(\psi(n)) < b_3 \log_2^2{x} \) hold for sufficiently large \( x \), where \( b_3 \) can be chosen to be any constant larger than \( 4e^3 \).

Thus, if \( n \in F \), then 
\( \omega(\psi(n)) \) is  \( O(\log_2^2{x}) \). To maximize $ h_{\psi}(n) $, we construct a prime sequence $ p_1 < p_2 < \cdots < p_{g(x)} $, where $ g(x) = (\log_2^2 x) $ and $ p_1 $ is the first prime greater than $ \log_2 x $. 
Obviously, for large \( x \), we have 
\[
\max {h_\psi(n)=\sum_{\substack{p| \psi \left ( n \right )  \\ p> \log_{2}x  }}\frac{1}{p}}\leq \sum_{i=1}^{g(x)} \frac{1}{p_i}  \leq \sum_{\log_2{x} < p < \log_2^3{x}} \frac{1}{p} \ll 1.
\]
Similarly
\[\max {h_\phi(n)}\ll 1.\]
Thus, by defining \( \phi(n) = n_1 \cdot n_2 \) and \( \psi(m) = m_1 \cdot m_2 \) as before, we derive that for \( n \in F \),
\begin{equation}
\begin{aligned}
\frac{\psi(\phi(n))}{\phi(n)} &=
\prod_{p \mid n_1} \left( 1 + \frac{1}{p}  \right)
\prod_{p \mid n_2} \left( 1 + \frac{1}{p} \right)\\
&\leq \prod_{p \leq \log_2{x}} \left( 1 + \frac{1}{p } \right) \exp(O(h_\phi(n))) \ll \log_3{x}\nonumber,
\end{aligned}
\end{equation}
and
\begin{equation}
\begin{aligned}
\frac{\phi(\psi(n))}{\psi(n)} &=
\prod_{p \mid m_1} \left( 1 - \frac{1}{p} \right) 
\prod_{p \mid m_2} \left( 1 - \frac{1}{p} \right)\\
&\geq \prod_{p \leq \log_2{x}} \left( 1 - \frac{1}{p} \right) \exp(-h_\psi(n)) \gg \frac{1}{\log_3{x}}.
\end{aligned}\nonumber
\end{equation}
Hence, \( I(n) \ll \log_3^2{x} \) for all \( n \in F \). From equation \eqref{eq:33}, we know that the set \(C_\phi \cup C_\psi \) has cardinality  \( O(\frac{x}{\sqrt{\log_3{x}}}) \), from which it follows that
\[
\sum_{n \in F} I(n) \leq \max_{n \in F} \{I(n)\} \cdot \#(C_\phi \cup C_\psi) \ll x \log_3^{3/2}{x}.
\]
Gathering together with \eqref{eq:21}, \eqref{eq:29} and \eqref{eq:35}, it follows that

\begin{equation}
\sum_{n \leq x} I(n) \leq \frac{6}{\pi^2}e^{2\gamma} c_0 x \log_3^2{x} + O(x \log_3^{3/2}{x}).\label{eq:36}\tag{4.13}
\end{equation}
By combining equations \eqref{eq:20} and\eqref{eq:36}, we complete the proof of equation \eqref{300} in Theorem 1.2.

\subsection{Proof of (\ref{301})}
We adopt the second method of proving Theorem 1 in \cite{ref8}, replacing the original ratio  
$\frac{\psi(n)}{n}$ with $\frac{\phi(n)}{\psi(n)}$   
  in the analysis. In view of \eqref{eq:18} and
\eqref{eq:28}, it follows that 

\[I(n)\geq \left(1 + o\left(1\right)\right) \frac{6}{\pi^2}  e^{2\gamma} \log_3^2 n \frac{\phi(n)}{\psi(n)}
\]
and
\[I(n) \le  \left(1 + o\left(1\right)\right)  \frac{6}{\pi^2}e^{2\gamma} \log_{3}^{2} n\frac{\phi(n)}{\psi(n)}
\]
 hold on a set of density 1. Therefore, it follows that on a set of density 1,

\[I(n) = \left(1 + o\left(1\right)\right) \frac{6}{\pi^2}e^{2\gamma} \log_3^{2} n\frac{\phi(n)}{\psi(n)}
=\left(1 + o\left(1\right)\right) \frac{6}{\pi^2}e^{2\gamma} \log_3^{2} n\prod_{p|n} \frac{p-1}{p+1} .\]
The proof of Theorem 1.2 is complete.
\section{Proof of Theorem 1.3 }
\subsection{Proof of (\ref{201})}

Similarly to the proof of Theorem 1.2, we can obtain
\begin{equation}
K(n) \geq 
\frac{6}{\pi^2}  e^{2\gamma} \log_3^2 x 
 + O\left( \log_3 x\cdot\log_4 x \right)
\label{eq:42}\tag{5.1}
\end{equation}
for $n \in \mathcal{A}.$
Additionally,
\begin{equation}
\sum_{n \leq x} K(n) \geq 
\sum_{\substack{n \in \mathcal{A} \\ n \leq x}}K(n) \geq \left(\frac{6}{\pi^2}  e^{2\gamma} \log_3^2 x 
 + O\left( \log_3 x\cdot\log_4 x \right)\right)\sum_{\substack{n \in \mathcal{A} \\ n \leq x}} 1.
\label{eq:43}\nonumber
\end{equation}

\begin{equation}
\sum_{n \leq x} K(n) \geq  \frac{6}{\pi^2}e^{2 \gamma} x \log _{3}^{2} x+O\left(x \log _{3} x\cdot\log_4 x\right) .
\label{eq:44}\tag{5.2}
\end{equation}
Following an approach analogous to $I(n)$, the interval partitioning method allows us to derive
\begin{equation}
K(n) \le  \frac{6}{\pi^2}e^{2\gamma} \log_3^2 x
\left(1 + O\left(\frac{1}{\sqrt{\log_3 x}}\right)\right)
\quad (n \in \mathcal{B}).\label{eq:45}\tag{5.3}
\end{equation}
Consequently,
\begin{equation}
\sum_{n \leq x} K(n) \leq \frac{6}{\pi^2}e^{2\gamma}  x \log_3^2{x} + O(x \log_3^{3/2}{x}).\label{eq:46}\tag{5.4}
\end{equation}
Combining \eqref{eq:46} and\eqref{eq:44} completes the proof of (\ref{201}) in Theorem 1.3.

\subsection{Proof of (\ref{302})}

Similarly to the proof of Theorem 1.2 and combining \eqref{eq:42} and \eqref{eq:45}, we obtain
\[K(n)\geq \left(1 + o\left(1\right)\right) \frac{6}{\pi^2}  e^{2\gamma} \log_3^2 n 
\]
and
\[K(n)\leq \left(1 + o\left(1\right)\right) \frac{6}{\pi^2}  e^{2\gamma} \log_3^2 n
\]
 hold on a set of density 1. Therefore, it follows that on a set of density 1,
\[K(n)= \left(1 + o\left(1\right)\right) \frac{6}{\pi^2}  e^{2\gamma} \log_3^2 n 
.\]
The proof of Theorem 1.3 is complete.

\end{document}